\documentclass{ifacconf}

\usepackage{graphicx}      
\usepackage{natbib}        
\usepackage{amsmath}
\usepackage{standalone} 
\usepackage{tikz}
\usepackage{makecell} 
\usepackage{booktabs}   
\usepackage{multirow}   

\usepackage[hyphens]{url}

\begin{document}
\begin{frontmatter}

\title{Mixed-integer programming formulations for optimal reconfiguration of supply chains }

\author[First,Second]{Lavinia M.P. Ghilardi} 
\author[Third]{Olga Walz} 
\author[Third]{Steffen Klosterhalfen} 
\author[First,Second]{Calvin Tsay}

\address[First]{Department of Computing, Imperial College London, London, SW7 2AZ, UK (e-mail: c.tsay@imperial.ac.uk, l.ghilardi@imperial.ac.uk)}
\address[Second]{Centre for Process Systems Engineering, Imperial College London, London, SW7 2AZ, UK}
\address[Third]{BASF, Ludwigshafen, 67056, Germany (e-mail: olga.walz@basf.com, steffen.klosterhalfen@basf.com)}

\begin{abstract}                %
Supply chains are interconnected networks
of processes and operations producing and delivering high-value products. These chains are increasingly subjected to structural changes from the energy transition and other external factors. To address this, this work develops mixed-integer programming formulations to identify optimal reconfigurations that preserve industrial operations and profitability.  We propose products and spatial neighborhoods to restrict the feasible search space and enable fast heuristic solutions. Furthermore, this restriction combines structural and product-based information, thus allowing to explore and define multiple reconfiguration scenarios. We demonstrate the approach using an agricultural waste case study, showing its ability to quickly produce good quality solutions.
\end{abstract}


\end{frontmatter}

\section{Introduction}

Supply chains are highly interconnected networks of processes and operations, spanning the journey and transformation of goods from raw materials to high-value products. Ensuring efficient operation of these supply chains is critical to guarantee the delivery of the products worldwide. Nevertheless, these supply chains are increasingly subjected to structural changes from external factors, such as the energy transition, which promotes sustainable product/process alternatives, or geopolitical events that reshape trade routes and resource availability. Identifying how to best reconfigure supply chains to respond to these structural changes is challenging, given their scale and interconnectedness.  

Current research in supply chain optimization has mainly focused on two types of problems: complete redesign of supply chains or short-term disruptions. The design problem has been extensively investigated over the years \citep{garciaYou2015}, focusing on resilient \citep{RibeiroBarbosa-Povoa2018, Arora2018} and agile supply chains \citep{agile, SHEKARIAN2020107438}. However, structural changes can affect existing supply chains, making complete redesigns impractical, and often infeasible. Another research direction addresses short-term disruptions \citep{Ivanov18102017, ORMS} and reactive responses \citep{BadejoIerapetritou2022}, but these fail to capture broader, planned structural changes such as the progressive transition toward sustainable solutions. A more realistic approach is to identify partial reconfigurations of supply chains that maintain industrial operations.  

Numerous computational decision-making tools have been proposed for supply chains, mainly focusing on routine operational decisions such as inventory management~\citep{burtea2024constrained,perez2021algorithmic}. 
On the other hand, to support decision-making during structural changes, alternative optimization tools are needed to help operators reconfigure some necessary part(s) of the network, maintaining operability and efficiency. Related \textit{reconfiguration problems} have been studied in other sectors, such as power transmission networks, reviewed by 
\citet{review_el}. 

Unlike the above, supply chains often involve heterogeneous technologies and multi-product transformations, making associated reconfiguration problems more specific and challenging. Multiple inventory and flow balances constraints must be satisfied, along with technology and capacity constraints. Mixed-Integer Linear Programming (MILP) has been widely applied to model and optimize supply chains \citep{LargeScaleSC, ovalle2024optimalreactiveoperationgeneral}, but remains computationally demanding for large-scale multi-product systems. To address computational scalability, approximation methods \citep{LargeScaleSC} and graph decompositions are usually used in complete design problems, while rolling horizon \citep{BadejoIerapetritou2022} and re-optimization approaches~\citet{review_el} and are used to tackle networks disruptions.

Given the above, this work focuses on optimizing (partial) reconfiguration of supply chains while maintaining operability and profitability. We first define a reconfiguration task as identifying the (multi-product) flows to be adjusted after a structural change. The problem is formulated as a MILP, with user-defined hyperparameters to guide desirable solutions in practice. Additionally, we propose several restrictions on the feasible reconnections based on neighborhood definitions. Beyond standard supply chain metrics based on spatial proximity, we introduce integrated neighborhoods, combining structural and product-based information to limit the search space and enable faster solution retrieval. Our contribution advances the current state-of-the-art by proposing efficient re-optimization strategies that exploit structural and product information.
The proposed methods are tested on instances from an agricultural case study \citep{LargeScaleSC}, showing that the neighborhood
definition allows to explore multiple reconfiguration scenarios efficiently.
\\
\section{Supply Chain Reconfiguration}

Supply chains are interconnected networks of processes and products flows, described by flow conservation equations, process and capacity constraints.  Reconfiguration of such networks may become necessary, for example, after the removal of a technology for environmental or other reasons. In such cases, the flows of goods must be adjusted to maintain operability and economic profitability. 

To model reconfiguration problems, we represent the supply chain as a \emph{spatial} graph $\mathcal{G}^S = (\mathcal{N}, \mathcal{A})$, following  \citet{LargeScaleSC}. Each node $i \in \mathcal{N}$ denotes a geographical site, and each directed arc $(i,j) \in \mathcal{A}$ 
represents a (multi-product) transport connection between sites. The nodes $i$ may host supply, demand, or technology units (or combinations thereof) for multiple products. Based on this representation, we first introduce the constraints describing the nominal network operation, followed by the formulation of the reconfiguration problem.

\subsection{Base Network Operation} \label{sec:basenet}

Given the possibility of having multiple products across the supply chain, each product $m \in \mathcal{M}$ is associated with a subset of nodes $\mathcal{N}_m$ and arcs $\mathcal{A}_m$ where it is present (i.e., where there is nonzero flow in nominal operation setting). 

Each node $i$ can host multiple suppliers $q \in \mathcal{S}_{i,m}$ and customers $c \in \mathcal{D}_{i,m}$, 
characterized by nonnegative flow quantities $s_q \ge 0$ and $d_c \ge 0$, respectively. 
Several technologies (processes) $t \in \mathcal{T}_i$ can be installed at each node. 
Each technology $t$ is associated with flow variables $x_{i,m,t}$ for each product $m$, where $x_{i,m,t} \ge 0$ indicates production and $x_{i,m,t} \le 0$ denotes consumption. 
Transport flows between nodes are represented by $f_{i,j,m} \ge 0$, indicating the transported quantity of $m$ along arc $(i,j) \in \mathcal{A}_m$. 
Accordingly, the flow conservation constraint for each product $m$ at each node $i$ is formulated as:
\begin{equation} 
\begin{aligned} 
&\sum_{q \in \mathcal{S}_{i,m}} s_{q} 
+ \sum_{t \in \mathcal{T}_i} x_{i,m,t} 
+ \sum_{j:(j,i)\in\mathcal{A}_m} f_{j,i,m} \\
&= \sum_{j:(i,j)\in\mathcal{A}_m} f_{i,j,m} 
+ \sum_{c \in \mathcal{D}_{i,m}} d_{c}, 
\quad \forall m \in \mathcal{M}, \ \forall i \in \mathcal{N}_m. 
\label{eq:mass_balance} 
\end{aligned} 
\end{equation}

Each technology $t$ operates according to a fixed input/output ratios $\hat{w}_{t,m}$, defined with respect to 
a reference $m_{\mathrm{ref},t}$ representing its main input product. 
Thus, technology flows satisfy:
\begin{equation}
\begin{split}
x_{i,m,t} &= \hat{w}_{t,m} \, x_{i,m_{\mathrm{ref},t},t}, \\
&\forall m \in \mathcal{M}, \ 
  \forall i \in \mathcal{N}_m, \ 
  \forall t \in \mathcal{T}_i.
\end{split}
\label{eq:tech_yield}
\end{equation}

Additionally, technology production is limited by its capacity $\hat{C}^{\text{tech}}_{i,t}$, defined with respect to $m_{\mathrm{ref},t}$:
\begin{equation} 
- x_{i,m_{\mathrm{ref},t},t} \leq \hat{C}^{\text{tech}}_{i,t}, 
\quad \forall i \in \mathcal{N}, \ \forall t \in \mathcal{T}_i.
\label{cap} 
\end{equation}
Similarly, the overall flows across the transport link are bounded by a fixed total capacity $\hat{f}^{\text{max}}_{i,j}$ for each arc $(i,j)$.
\begin{equation} 
\sum_{m \in\mathcal{M}} f_{j,i,m} \leq \hat{f}^{\text{max}}_{i,j}, 
\quad \forall (i,j) \in \mathcal{A}.
\label{cap} 
\end{equation}

The optimal base network configuration is typically found by solving a combinatorial optimization problem with constraints \eqref{eq:mass_balance}--\eqref{cap} and an economic objective function. 
Nevertheless, in this work we define a general reconfiguration problem that only requires a \textit{feasible} supply chain as a starting point. 
Therefore, we do not require the initial supply chain (before disruption) to be optimal, but assume it is feasible with respect to constraints \eqref{eq:mass_balance}--\eqref{cap}. 

\subsection{Reconfiguration Problem} \label{sec:rec}

This section presents the Mixed-Integer Linear Programming (MILP) formulation of the reconfiguration problem. Our presentation focuses on the reconfiguration following the removal of a technology, 
where products flows must be adjusted to maintain operability and economic performance. 
The proposed formulation, however, is generic to other reconfiguration scenarios involving node removals, 
such as the deactivation of suppliers or customers. 

In graph terms, the reconfiguration problem corresponds to the removal of a technology $r$ from a node $i_r$, which is enforced in the formulation by fixing the associated technology flows to zero:
\begin{equation}
x_{i_r,m,r} = 0
\quad \forall m \in \mathcal{M}_{i_r}.
\label{eq:tech_removed}
\end{equation}
Following this definition, all other entities associated with node $i_r$, such as supply, demand, or other technologies $t \neq r$, are preserved.  
In other words, node $i_r$ remains active, but loses the production of specific products previously generated by technology $r$.

The reconfiguration task is then to identify new potential transport connections, represented by binary variables $z^{\text{rec}}_{i,j}$. These (re)connections can be selected among the candidate arcs $\mathcal{A}^{\text{cand}}$. This set ensures product consistency between nodes $i$ and $j$, and can be defined/restricted via the neighborhoods presented later in Section~\ref{sec:chem_neig}. 

Product flows are allowed for an arc only if the corresponding reconnection is selected:
\begin{equation}
\sum_{m \in \mathcal{M}_{i,j}} f_{i,j,m} \le 
z^{\text{rec}}_{i,j} \, \hat{f}^{\text{max}}_{i,j},
\quad \forall (i,j) \in \mathcal{A}^{\text{cand}}.
\label{eq:rec_arc_flow}
\end{equation}
Accordingly, the flow conservation constraint \eqref{eq:mass_balance} 
is reformulated in \eqref{eq:mass_balance_rec} to account for the newly introduced reconnections. For each product $m$, the set $\mathcal{A}_m^{\text{cand}}$ is defined as the subset of candidate arcs $(i,j) \in \mathcal{A}^{\text{cand}}$ where product $m$ is present at both nodes (i.e., $i \in \mathcal{N}_m \land j \in \mathcal{N}_m$).

\begin{equation} 
\begin{aligned} 
&\sum_{q \in \mathcal{S}_{i,m}} s_{q} 
+ \sum_{t \in \mathcal{T}_i} x_{i,m,t} 
+ \sum_{j:(j,i)\in\mathcal{A}_m \cup \mathcal{A}_m^{\text{cand}}} f_{j,i,m} \\
&= \sum_{j:(i,j)\in\mathcal{A}_m \cup \mathcal{A}_m^{\text{cand}}} f_{i,j,m} 
+ \sum_{c \in \mathcal{D}_{i,m}} d_{c} 
\quad \forall m \in \mathcal{M}, \ \forall i \in \mathcal{N}_m
\label{eq:mass_balance_rec} 
\end{aligned} 
\end{equation}

To ensure practical relevance, we now introduce additional parameters that allow reconfiguration preferences to be specified. In particular, the maximum number of reconnections $\hat{N}^{\text{rec}}$ and the penalty parameter $\hat{\phi}^{\text{d}}$ for deviations from baseline customer demand can be specified. 
Constraint~\eqref{eq:rec_arc_limit} limits the number of new transport links according to $\hat{N}^{\text{rec}}$:
\begin{equation}
\sum_{(i,j)\in\mathcal{A}^{\text{cand}}} z^{\text{rec}}_{i,j} \le \hat{N}^{\text{rec}}.
\label{eq:rec_arc_limit}
\end{equation}

To penalize deviations from the baseline customer demand $\hat{d}^{\text{base}}_{c}$, 
a nonnegative auxiliary variable $\delta^{\text{d}}_{c}$ is introduced:
\begin{equation}
\delta^{\text{d}}_{c} \ge d_{c} - \hat{d}^{\text{base}}_{c}, \quad 
\delta^{\text{d}}_{c} \ge -\big(d_{c} - \hat{d}^{\text{base}}_{c}\big),
\quad \forall c \in \mathcal{D}.
\label{eq:delta_dem}
\end{equation}

The objective is to maximize overall economic profit, while penalizing deviations from the baseline demand through the weighting parameter $\hat{\phi}^{\text{d}}$:
\begin{align}
\max \quad &
\sum_{c \in \mathcal{D}} \hat{c}^{\text{d}}_c \, d_{c}
- \sum_{q \in \mathcal{S}} \hat{c}^{\text{s}}_q \, s_{q}
- \sum_{i\in\mathcal{N}} \sum_{t \in\mathcal{T}_i}
  \hat{c}^{\text{tech}}_t \, (-x_{i,m_{\mathrm{ref},t},t})
\nonumber \\ &
- \sum_{m\in\mathcal{M}} \sum_{(i,j)\in\mathcal{A}_m \cup \mathcal{A}_m^{\text{cand}}}
  \hat{c}^{\text{trans}}_{i,j,m} \, f_{i,j,m}
- \hat{\phi}^{\text{d}} \sum_{c\in\mathcal{D}} \delta^{\text{d}}_c,
\label{eq:obj_reconfig}
\end{align}
where $\hat{c}^{\text{d}}_c$ denotes the unit revenue from customer demand, 
$\hat{c}^{\text{s}}_q$ the supply cost, 
$\hat{c}^{\text{tech}}_t$ the operating cost of technology $t$, 
and $\hat{c}^{\text{trans}}_{i,j,m}$ the transport cost of material $m$ along arc $(i,j)$. 
Note that the objective function \eqref{eq:obj_reconfig} without the final penalty term can be selected to be consistent with the objective of the original base case supply chain design. 

Overall, the reconfiguration problem determines the set of new transport links and corresponding material flows that maximize the total profit, 
while maintaining customer demand in~\eqref{eq:obj_reconfig} and satisfying all operational constraints 
\eqref{eq:tech_yield}--\eqref{eq:delta_dem}, including flow conservation constraints, technology capacities, and reconfiguration limits. 
While this initial work focuses on reconfiguration of a single structural disturbance, the framework may be generalized to establish reconfigurations for a sequence of disturbances.

\section{Neighborhood Definition}
\label{sec:chem_neig}

The reconfiguration problem identifies reconnections, i.e., new product transport arcs to preserve operations and profitability. Potential arcs $z^{\text{rec}}_{i,j}$ are chosen from a candidate set $\mathcal{A}^{\text{cand}}$, restricted to ensure material consistency and avoid existing links $\mathcal{A}$. Without restrictions, this set equals the full combinatorial superstructure:
\begin{equation}
\begin{aligned}
\mathcal{A}^{\text{full}}
= \big\{ (i,j) :\;&
i, j \in \mathcal{N},\ j \neq i,\ (i,j) \notin \mathcal{A},\\
&\exists m \in \mathcal{M}:\ i \in \mathcal{N}_m \land j \in \mathcal{N}_m
\big\},
\end{aligned}
\label{eq:neigh_full}
\end{equation}
where $\mathcal{A}^{\text{full}}$ can grow quadratically with the
number of sites. To both reduce the MILP search space and allow specification of the degree of reconfiguration, we restrict reconnections to a node neighborhood $\mathcal{N}_k(i_r,r)$ around the removed technology $r$ at node $i_r$ within a distance $k$:
\begin{equation}
\mathcal{A}^{\text{cand}}
= \{ (i,j) \in \mathcal{A}^{\text{full}} :
i,j \in \mathcal{N}_k(i_r, r) \}.
\label{eq:neigh_struct}
\end{equation}
\\
Different definitions of $\mathcal{N}_k(i_r,r)$ are presented in the following to balance computational cost and solution quality.
Three neighborhood types are defined to progressively incorporate different information of the supply chain: (i) \emph{spatial}, (ii) \emph{product}, and (iii) integrated \emph{ spatial-product}. 

\begin{figure}[t]
  \centering
\includegraphics[width=0.8\columnwidth]{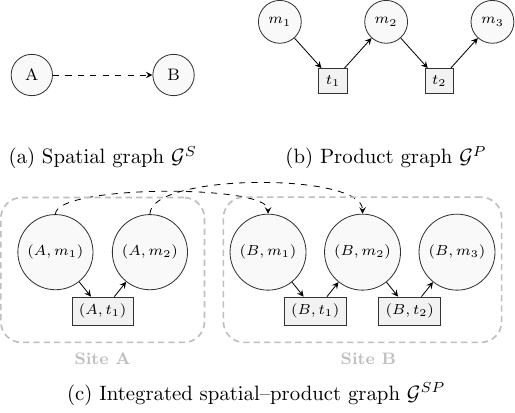}
  \caption{Supply chain graph representations used for neighborhood definitions}

  \label{fig:graphSC}
\end{figure}

To build these neighborhoods, the supply chain can be viewed through three complementary graphs, shown in Fig.~\ref{fig:graphSC}.  
The \emph{spatial}  graph $\mathcal{G}^S$ describes existing transport links between geographical sites.  
The \emph{product} graph $\mathcal{G}^P$ connects products and processes.  
Finally, the integrated \emph{ spatial-product} graph $\mathcal{G}^{SP}$ combines both perspectives, capturing how products can be transformed and transported across the supply chain.

\subsection{Spatial Neighborhoods}

Spatial neighborhoods rely on the \emph{spatial} graph 
$\mathcal{G}^S = (\mathcal{N}, \mathcal{A})$, 
where nodes represent geographical sites, and arcs denote existing transport links.  
Let $\Delta^{\text{geo}}$ be a specified geographical threshold.  
Then, the $\Delta$-geographical neighborhood around site $i_r$ is defined as
\begin{equation}
\mathcal{N}^{S}_\Delta(i_r)
= \{ j \in \mathcal{N} : D(i_r,j) \le \Delta^{\text{geo}} \}.
\label{eq:neigh_geo}
\end{equation}

Alternatively, the graph connectivity can be exploited through the shortest-path distance 
$d^S(\cdot,\cdot)$ between sites $i$ and $j$ in the undirected version of $\mathcal{G}^S$.  
The corresponding $k$-step spatial neighborhood is $\mathcal{N}^{S,\text{step}}_k(i_r)
= \{ j \in \mathcal{N} : d^S(i_r,j) \le k \}$.
Both definitions rely solely on the spatial structure of the supply chain as defined by distances, and therefore disregard transformations and value additions between materials.

\subsection{Product Neighborhoods}

Spatial neighborhoods neglect the heterogeneous, multi-product nature of supply chains.  
To address this, we define a \emph{product} neighborhood that captures how products and technologies are interconnected through the supply chain.  
Intuitively, this definition identifies the products associated with the removed technology through a `distance' in number of transformation steps. This enables heuristically restricting reconfigurations in terms of product proximity, rather than physical distance.

Formally, the product neighborhood is constructed on the technology-material connectivity of the \emph{product} graph $\mathcal{G}^P = (\mathcal{V}^P,\mathcal{E}^P)$. 
The bipartite graph $\mathcal{G}^P$, illustrated in Fig.~\ref{fig:graphSC}, contains products $m \in \mathcal{M}$ and technologies $t \in \mathcal{T}$ as nodes, with arcs $(m,t)$ if $m$ is consumed by $t$, and $(t,m)$ if $m$ is produced. 

The shortest-path distance $d^P(\cdot,\cdot)$ on the undirected version of $\mathcal{G}^P$ defines the set of products reachable from the removed technology $r$ within $k$ processing steps:
\begin{equation}
\mathcal{M}^{P,\text{reach}}_{k}(r)
= \{ m \in \mathcal{M} : d^P(r,m) \le k \}.
\label{eq:chem_n_mat}
\end{equation} 

A product that can be formed within $\le k$ process transformations from technology $r$ is considered \emph{product reachable}. 

Using this definition, we specify the product neighborhood relevant for the reconfiguration as 
sites $j \in \mathcal{N}$ compatible with reachable products $\mathcal{M}^{P,\text{reach}}_{k}(r)$:
\begin{equation}
\mathcal{N}^{P}_k(r)
= \{ j \in \mathcal{N}: \exists m \in \mathcal{M}_j \cap \mathcal{M}^{P,\text{reach}}_{k}(r)\}.
\label{eq:chem_n}
\end{equation}

This formulation effectively maps product-technology connectivity back onto the spatial network, identifying sites that can provide compatible products for reconfiguration.

\subsection{Spatial-Product Neighborhoods}

Product neighborhoods describe transformation connectivity, but disregard transport feasibility between sites.  
To capture both aspects simultaneously, we introduce an integrated \emph{spatial-product}  neighborhood.  
Conceptually, this definition integrates measures of how materials and processes are connected when both process transformations and transport links are considered.  
This combined view enables realistic reconfiguration decisions based on both process and transport pathways.

Formally, the integrated neighborhood is defined on the 
\emph{spatial-product graph} 
$\mathcal{G}^{SP} = (\mathcal{V}^{SP}, \mathcal{E}^{SP})$, 
shown in Fig.~\ref{fig:graphSC}.  
It provides a representation similar to the spatial superstructure of \cite{SHAO2023108102}.  

The node set $\mathcal{V}^{SP}$  includes site-material pairs $(i,m)$ and site-process pairs $(i,t)$. 
The arcs $\mathcal{E}^{SP}$ represent either (i) transport of materials between sites $\mathcal{E}^{SP,\text{tr}}$ 
or (ii) transformations within processes $\mathcal{E}^{SP,\text{in}}$ and $\mathcal{E}^{SP,\text{out}}$. 
\begin{align}
\mathcal{E}^{SP,\mathrm{tr}} &= \{ ((i,m),(j,m)) : (i,j)\in\mathcal{A}_m \},\\
\mathcal{E}^{SP,\mathrm{in}} &= \{ ((i,m),(i,t)) : i\in\mathcal{N},\ t\in\mathcal{T}_i, \notag\\
&\phantom{= \{ ((i,m),(i,t)) :} m\in\mathcal{M}_t^{\mathrm{in}} \},\\
\mathcal{E}^{SP,\mathrm{out}} &= \{ ((i,t),(i,m)) : i\in\mathcal{N},\ t\in\mathcal{T}_i, \notag\\
&\phantom{= \{ ((i,t),(i,m)) :} m\in\mathcal{M}_t^{\mathrm{out}} \},
\end{align}
where ${M}_t^{\mathrm{in}}$ represents the sets of the input products to technology $t$ and ${M}_t^{\mathrm{out}}$ the corresponding outlet products.
This representation preserves the bipartite structure between products $m$ and technologies $t$.
A single step on $\mathcal{G}^{SP}$ therefore corresponds to a product transport or a technology conversion, 
and the shortest-path distance $d^{SP}(\cdot,\cdot)$ defines a unified measure of connectivity.

Given the removed technology $r$ at site $i_r$, the set of reachable products within $k$ spatial-product steps can be computed similarly to \eqref{eq:chem_n_mat} using $d^{SP}, (\cdot,\cdot)$, while the projection of these products onto a set of compatible sites can be computed  \eqref{eq:chem_n}. 
However, this sequential projection effectively first applies the product logic 
and then the spatial mapping.

In contrast, we directly identify the neighborhood by using
spatial and product connectivity simultaneously:
\begin{equation}
\mathcal{N}^{SP}_k(i_r,r)
=\!\left\{
j \in \mathcal{N} :
\begin{array}{l}
\exists\, \ m \in \mathcal{M}, \\[2pt]
d^{SP}((i_r,r),(j,m)) \le k
\end{array}
\!\right\}.
\label{eq:chem_nsc}
\end{equation}

This definition leverages the full structure of $\mathcal{G}^{SP}$, 
capturing transport and product transformation jointly.

In summary, spatial neighborhoods constrain reconfigurations by 
transport logistics, product neighborhoods by product-technology connectivity, 
and integrated spatial-product neighborhoods by both.  
These neighborhoods restrict the reconfiguration search space, 
improving tractability by exploiting the structure of supply chains.

\section{Results and Discussion}

The proposed formulations are applied to an agricultural waste management case study \citep{SAMPAT2019352}, where the supply chain transforms agricultural residues into value-added products. The products are transported across croplands, dairy and beef farms in Wisconsin, along with external consumers. Instances are derived from the network design problem of \cite{LargeScaleSC}, including 9 products (3 waste-derived raw materials and 6 final products) and 3 possible transformation processes.

Our analysis considers the reconfiguration following the removal of a process $r$ at site $i_r$. Different neighborhood definitions are compared in terms of computational effort, solution quality, and operational disruption (e.g., deviation of customer demand). The maximum number of reconnections $\hat{N}^{\text{rec}}$ is varied to provide a set of alternative solutions for industrial decision-making. MILPs are implemented in Pyomo \citep{bynum2021pyomo} and solved using Gurobi 12.0.2 \citep{gurobi} on an Apple M4 Pro CPU with a MIP gap of 0.1\% and a time limit of 3600 seconds.

We consider four reconfiguration scenarios: (i) \textbf{Full (F)}, without restrictions on $\mathcal{A}^{\text{cand}}$; (ii) \textbf{Spatial (N-S)}, limited by geographical distance $\Delta^{\text{geo}}$; (iii) \textbf{Product (N-P)}, based on reachable products within $k$- transformations; and (iv) \textbf{Integrated Spatial-Product (N-SP)}, considering $k$- technology and transports steps across the supply chain in a joint manner. The corresponding candidate sets $\mathcal{A}^{\text{cand}}$ in each reconfiguration problem are derived from equations \eqref{eq:neigh_full}, \eqref{eq:neigh_struct}, \eqref{eq:chem_n}, \eqref{eq:chem_nsc}, respectively. 
The parameter $\hat{\phi}^{\text{d}}$ in \eqref{eq:obj_reconfig} balances economic profitability and degree of reconfiguration. We set it empirically to twice the maximum demand revenue to discourage demand deviation, though its value can be tuned based on user preference.
\subsection{Case Study 1: Small-Scale Supply Chain}
\begin{figure}[t]
  \centering
 \includegraphics[width=0.4\textwidth]{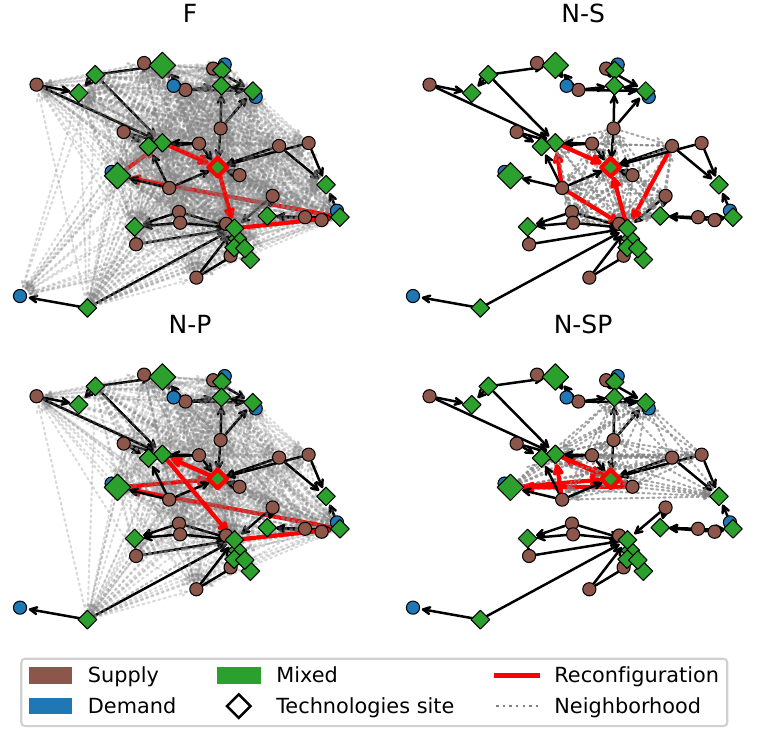}
  \caption{Small case study: reconfigurations for different neighborhoods ($\hat{N}^{rec}=5$)}
  \label{fig:rec_VC}
\end{figure}

We first consider a supply chain of 46 nodes, with the spatial structure shown in Fig.~\ref{fig:rec_VC}. As a naive baseline, simply removing technology $r$ without any reconfiguration of the supply chain yields an objective function value of 16.40, whereas allowing reconfiguration improves this value  to 17.54, as described below.
Dashed gray lines depict the candidate reconnections $\mathcal{A}^{\text{cand}}$ defined by each neighborhood, where we empirically observe that each neighborhood definition produces a significantly different set of candidate arcs for reconfiguration. Parameters are set to $k=1$ for N-P (direct input/output products from the removed process) and $k=2$ for N-SP (one product and one transport step). We empirically tune the distance-based, geographical threshold $\Delta^{\text{geo}}$ in N-S to yield candidate sets of comparable sizes with respect to N-SP to assess the impact of integrating product-based information within the spatial setting. These parameters were selected to restrict the solution space and improve computational time while enabling a sufficient degree of reconfiguration. In practice, they can be tuned adaptively by progressively increasing $k$ and warm-starting the solution.

Results for the different scenarios are summarized in Table \ref{tab:reconfig_neigh}, showing the cardinality of \textbf{$|\mathcal{A}^{\text{cand}}|$}, the customer demand deviation after the reconfiguration, the objective function and the problem solution time. 
The reduction in the size of \textbf{$|\mathcal{A}^{\text{cand}}|$} induced by N-P is less significant compared to other neighborhood definitions. This reflects the dense product distribution of this specific supply chain, where most materials are present at multiple sites. As a result, this restriction produces high-quality solutions: the objective function is similar to searching over $\mathcal{A}^{\text{full}}$ (denoted F), but with lower computational time. N-SP and N-S exhibit moderate demand and objective function deviations from F, but lower computational times. Between the two, N-SP yields better solutions.  
When the maximum number of reconfigurable arcs $\hat{N}^{\text{rec}}$ is limited to 2, variations in both the objective function and customer demand deviation are relatively small across the different neighborhood definitions. Increasing $\hat{N}^{\text{rec}}$ to 5 amplifies the impact of the neighborhood choice, with the best objective value equal to 17.54. However, we found increasing $\hat{N}^{\text{rec}}$ beyond 5 to not significantly improve solutions of the reconfiguration problem. The resulting optimal connections are shown in red in Fig.~\ref{fig:rec_VC} for $\hat{N}^{\text{rec}}=5$. 

\begin{table}[ht]
\centering
\caption{Small case study: comparison of reconfiguration across neighborhoods.}
\footnotesize
\setlength{\tabcolsep}{4pt}
\renewcommand{\arraystretch}{0.95}
\begin{tabular}{@{}c c c c c c@{}}
\toprule
\makecell[c]{\textbf{$\hat{N}^{\text{rec}}$}} &
\makecell[c]{\textbf{Neigh.}} &
\makecell[c]{\textbf{$|\mathcal{A}^{\text{cand}}|$}} &
\makecell[c]{\textbf{Demand}\\\textbf{Dev. [\%]}} &
\makecell[c]{\textbf{Obj.}\\\textbf{Fun.}} &
\makecell[c]{\textbf{Time}\\\textbf{[s]}} \\
\midrule
\multirow{4}{*}{2}
 & F     & 1478 & 2.5 & 16.83 & 138.62 \\
 & N--S  & 141  & 2.8 & 16.68 & 0.25 \\
 & N--P  & 941  & 2.5 & 16.83 & 52.38 \\
 & N--SP & 133  & 2.5 & 16.83 & 0.35 \\
\midrule
\multirow{4}{*}{5}
 & F     & 1478 & 1.0 & 17.53 & 2.63 \\
 & N--S  & 141  & 2.4 & 16.87 & 0.01 \\
 & N--P  & 941  & 1.0 & 17.53 & 1.15 \\
 & N--SP & 133  & 2.1 & 17.01 & 0.01 \\
\midrule
\multirow{4}{*}{10}
 & F     & 1478 & 1.0 & 17.54 & 0.06 \\
 & N--S  & 141  & 2.4 & 16.88 & 0.01 \\
 & N--P  & 941  & 1.0 & 17.54 & 0.04 \\
 & N--SP & 133  & 2.1 & 17.01 & 0.01 \\
\bottomrule
\end{tabular}
\label{tab:reconfig_neigh}
\end{table}
The results indicate that smaller candidate sets, such as those defined by the N-S and N-SP neighborhood, retain high quality solutions despite the reduced search space. This demonstrates that carefully chosen neighborhoods preserve essential demand operations and profitability in the reconfiguration, while decreasing the computational time. Overall, the comparison across neighborhoods highlights how combining spatial and product information in N-SP leads to the best compromise between the quality of the solution and the computational time.

\subsection{Case Study 2: Large-Scale Supply Chain}

\begin{figure}[ht]
  \centering
 \includegraphics[width=0.3\textwidth, trim=1cm 0pt 1cm 1cm, clip]{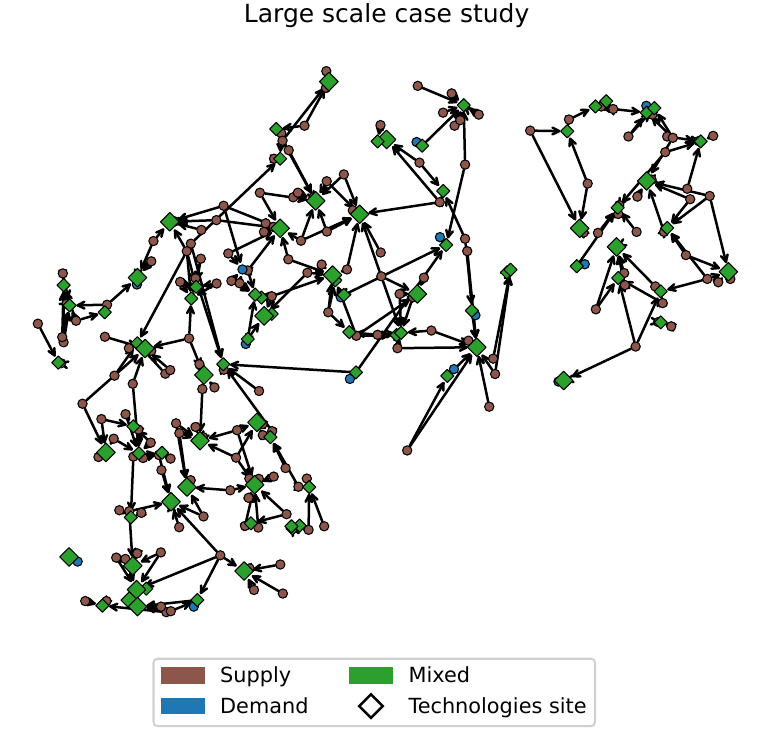}
  \caption{Large case study: supply chain graph ($\hat{N}^{rec}=5$)}
  \label{fig:large_VC}
\end{figure}

To assess scalability of the proposed methods, we now consider a larger instance with 278 nodes, with structure shown in Fig. \ref{fig:large_VC}. In this setting, naively removing the technology decreases the objective function to 157.63. 

Table~\ref{tab:large_reconfig_neigh} compares the results under varying $\hat{N}^{\text{rec}}$ and neighborhood definitions. Given the larger size of the supply chain, the set $|\mathcal{A}^{\text{cand}}|$ contains 43275 candidate arcs in the full reconfiguration scenario, which cannot be solved in the given time limit. The neighborhood restrictions N-S and N-SP substantially reduce the number of candidate arcs $|\mathcal{A}^{\text{cand}}|$ to only 858 and 854, respectively, enabling complete MILP solution. As in the small case study, the search space reduction provided by N-P is lower due to dense product distribution, with 35808 arcs. 
Though this reduction is moderate, the resulting MILP is solved in 52.38s, suggesting the neighborhood definition still produces a better-posed optimization formulation. 

When $\hat{N}^{rec}=5$, F and N-P achieve better objective function and lower demand deviation than N-S and N-SP. However, the computational time required for solving F and N-P is higher, due to the increased size of $|\mathcal{A}^{\text{cand}}|$. 
Therefore, we observe some tradeoff between the computational tractability enabled by a restriction, and the resulting profitability found. 
Overall, deviation from F in terms of the objective function is below 3\% in N-S and N-SP for all values of $\hat{N}^{rec}$, with N-SP offering better quality solutions than N-S. 

Given these results, the combined N-SP  definition appears to offer the best balance between tractability and optimality, even for large instances. Results depend on the studied system, where widespread goods across sites reduce the impact of product information. In high-value product chains, more intermediate steps yield sparser product distributions and stronger reductions in problem size. However, curtailed demand for valuable products can increase economic  losses, making reconfiguration critical.

\begin{table}[ht]
\centering
\caption{Large case study: comparison of reconfiguration across neighborhoods.}
\footnotesize
\setlength{\tabcolsep}{4pt}
\renewcommand{\arraystretch}{0.95}
\begin{tabular}{@{}c c c c c c@{}}
\toprule
\makecell[c]{\textbf{$\hat{N}^{\text{rec}}$}} &
\makecell[c]{\textbf{Neigh.}} &
\makecell[c]{\textbf{$|\mathcal{A}^{\text{cand}}|$}} &
\makecell[c]{\textbf{Demand}\\\textbf{Dev. [\%]}} &
\makecell[c]{\textbf{Obj.}\\\textbf{Fun.}} &
\makecell[c]{\textbf{Time}\\\textbf{[s]}} \\
\midrule
\multirow{4}{*}{2}
& F     & 43275 & 1.17 & 159.09 & 3600 \\
& N--S  & 858   & 1.17 & 159.06 & 3.50 \\
& N--P & 35808 & 1.15 & 159.08 & 52.38 \\
& N--SP  & 854   & 1.16 & 159.08 & 3.77 \\

\midrule
\multirow{4}{*}{5}
 & F     & 43275 & 0.19 & 163.09 & 3600 \\
 & N--S  & 858   & 0.91 & 160.17 & 0.26 \\
& N--P & 35808 & 0.19 & 163.09 & 2330.8 \\
 & N--SP  & 854   & 0.54 & 161.71 & 88.24 \\

\midrule
\multirow{4}{*}{10}
 & F     & 43275 & 0.00 & 164.03 & 506.25 \\
 & N--S  & 858   & 0.89 & 160.27 & 0.04 \\
  & N--P & 35808 & 0.00 & 164.01 & 118.11 \\
 & N--SP  & 854   & 0.41 & 162.30 & 0.23 \\

\bottomrule
\end{tabular}
\label{tab:large_reconfig_neigh}
\end{table}

\section{Conclusion}

Reconfigurations of supply chains are essential to adapt to sustainability transitions and other structural changes. This work develops MILP formulations for optimal reconfiguration and neighborhoods to restrict the feasible search space. Application to an agricultural waste case study shows that combining spatial and product-based information significantly reduces computational time, while maintaining high-quality solutions. 
This physically motivated restriction allows exploration of multiple reconfiguration scenarios even in larger case problems. Future work will extend the approach to diverse supply chains and explore scalable optimization and heuristic strategies. 

\begin{ack}
Support from BASF and the EPSRC (EP/X025292/1) is gratefully acknowledged. CT was supported by a BASF/Royal Academy of Engineering Senior Research Fellowship. 
The authors thank Dr Antonio del R\'io Chanona for the insightful feedback and suggestions during the development of this research.

\end{ack}

\bibliography{ifacconf}            
\end{document}